\documentstyle[txmac,a4,
amssymb,%
twoside,%
nocaphead,%
epsf,%
mypic,times,mathptm]{article}

\advance\oddsidemargin by -1.8cm
\advance\evensidemargin by -1.8cm
\advance\textwidth by 3.6cm

\def\mynewtheo#1#2{%
\newtheorem{@#1}{#2}[section]%
\newenvironment{#1}{\begin{@#1}\rm}{\end{@#1}}}

\mynewtheo{lemma}{Lemma}
\mynewtheo{exer}{Exercise}
\mynewtheo{theo}{Theorem}
\mynewtheo{remark}{Remark}
\mynewtheo{defi}{Definition}
\mynewtheo{conj}{Conjecture}
\mynewtheo{corr}{Corollary}
\mynewtheo{prop}{Proposition}
\mynewtheo{question}{Question}
\mynewtheo{exam}{Example}

\newenvironment{theorem}{\begin{theo}}{\end{theo}}

\newenvironment{eqn}{\begin{equation}}{\end{equation}\ignorespaces}

\begin{document}

\makeatletter

\def\bysame{\same[\kern2cm]\,}
\def\lra{\longrightarrow}
\def\Md{\max\deg}
\def\md{\min\deg}
\def\qed{\hfill\@mt{\Box}}
\def\@mt#1{\ifmmode#1\else$#1$\fi}
\def\qqed{\hfill\@mt{\Box\enspace\Box}}
\def\TM{$^\text{\raisebox{-0.2em}{${}^\text{TM}$}}$}

\let\ap\alpha
\let\tl\tilde
\let\lm\lambda
\let\Lm\Lambda
\let\dl\delta
\let\Dl\Delta
\let\eps\varepsilon
\let\gm\gamma
\let\bt\beta
\def\cF{{\cal F}}
\def\bZ{{\Bbb Z}}
\def\bQ{{\Bbb Q}}
\def\bN{{\Bbb N}}
\def\bR{{\Bbb R}}
\def\bC{{\Bbb C}}
\def\hK{{\hat K}}
\def\hKhK{\hK\#\hK}
\let\sg\sigma
\let\dt\det
\let\op\oplus

\def\epsfs#1#2{{\catcode`\_=11\relax\ifautoepsf\unitxsize#1\relax\else
\epsfxsize#1\relax\fi\epsffile{#2.eps}}}
\def\epsfsv#1#2{{\vcbox{\epsfs{#1}{#2}}}}
\def\vcbox#1{\setbox\@tempboxa=\hbox{#1}\parbox{\wd\@tempboxa}{\box
  \@tempboxa}}
\def\p{\epsfsv{2cm}}

\def\@test#1#2#3#4{%
  \let\@tempa\go@
  \@tempdima#1\relax\@tempdimb#3\@tempdima\relax\@tempdima#4\unitxsize\relax
  \ifdim \@tempdimb>\z@\relax
    \ifdim \@tempdimb<#2%
      \def\@tempa{\@test{#1}{#2}}%
    \fi
  \fi
  \@tempa
}

\def\go@#1\@end{}
\newdimen\unitxsize
\newif\ifautoepsf\autoepsftrue

\unitxsize4cm\relax
\def\epsfsize#1#2{\epsfxsize\relax\ifautoepsf
  {\@test{#1}{#2}{0.1 }{4   }
		{0.2 }{3   }
		{0.3 }{2   }
		{0.4 }{1.7 }
		{0.5 }{1.5 }
		{0.6 }{1.4 }
		{0.7 }{1.3 }
		{0.8 }{1.2 }
		{0.9 }{1.1 }
		{1.1 }{1.  }
		{1.2 }{0.9 }
		{1.4 }{0.8 }
		{1.6 }{0.75}
		{2.  }{0.7 }
		{2.25}{0.6 }
		{3   }{0.55}
		{5   }{0.5 }
		{10  }{0.33}
		{-1  }{0.25}\@end
		\ea}\ea\epsfxsize\the\@tempdima\relax
		\fi
		}

\let\old@tl\~
\def\~{\raisebox{-0.8ex}{\tt\old@tl{}}}

\author{A. Stoimenow\footnotemark[1]\\[2mm]
\small Research Institute for Mathematical Sciences, \\
\small Kyoto University, Kyoto 606-8502, Japan\\
\small e-mail: {\tt stoimeno@kurims.kyoto-u.ac.jp}\\
\small WWW: {\hbox{\web|http://www.kurims.kyoto-u.ac.jp/~stoimeno/|}}
}

{\def\thefootnote{\fnsymbol{footnote}}
\footnotetext[1]{Supported by 21st Century COE Program.}
\def\thefootnote{}
\footnotetext[2]{
\ry{1.3em}\em{AMS subject classification:} 57M25 (primary),
11E04, 11T30 (secondary), \\
\rx{1.8em}%
\em{Key words:} quadratic form, finite commutative group,
4-genus, signature, slice knot, Alexander polynomial.
}}

\title{\large\bf \uppercase{Some examples related to knot sliceness}
\\[4mm]
\phantom{\small\it This is a preprint. I would be grateful
for any comments and corrections!}}

\date{\phantom{\large Current version: \curv\ \ \ First version:
\makedate{13}{4}{2003}}}

\maketitle

\long\def\@makecaption#1#2{%
   \vskip 10pt
   {\let\label\@gobble
   \let\ignorespaces\@empty
   \xdef\@tempt{#2}%
   }%
   \ea\@ifempty\ea{\@tempt}{%
   \setbox\@tempboxa\hbox{%
      \fignr#1#2}%
      }{%
   \setbox\@tempboxa\hbox{%
      {\fignr#1:}\capt\ #2}%
      }%
   \ifdim \wd\@tempboxa >\captionwidth {%
      \rightskip=\@captionmargin\leftskip=\@captionmargin
      \unhbox\@tempboxa\par}%
   \else
      \hbox to\captionwidth{\hfil\box\@tempboxa\hfil}%
   \fi}%
\def\fignr{\small\sffamily\bfseries}%
\def\capt{\small\sffamily}%

\newdimen\@captionmargin\@captionmargin2cm\relax
\newdimen\captionwidth\captionwidth\hsize\relax

\let\reference\ref
\def\eqref#1{(\protect\ref{#1})}

\def\proof{\@ifnextchar[{\@proof}{\@proof[\unskip]}}
\def\@proof[#1]{\noindent{\bf Proof #1.}\enspace}

\def\myfrac#1#2{\raisebox{0.2em}{\small$#1$}\!/\!\raisebox{-0.2em}{\small$#2$}}
\def\abstractname{}

\def\br#1{\left\langle#1\right\rangle}

\@addtoreset {footnote}{page}

\renewcommand{\section}{%
   \@startsection
         {section}{1}{\z@}{-1.5ex \@plus -1ex \@minus -.2ex}%
               {1ex \@plus.2ex}{\large\bf}%
}
\renewcommand{\@seccntformat}[1]{\csname the#1\endcsname .
\quad}

{\let\@noitemerr\relax
\vskip-2.7em\kern0pt\begin{abstract}
\noindent{\bf Abstract.}\enspace
It is known that the linking form on the 2-cover of slice knots
has a metabolizer. We show that several weaker conditions, or
some other conditions related to sliceness, do not imply the
existence of a metabolizer. We then show how the Rudolph-Bennequin
inequality can be used indirectly to prove that some knots are
not slice. \\[2mm]
\end{abstract}
}

\section{Introduction and statement of results}

Algebraic topology has attempted the study of topological
equivalence problems by means of algebraic invariants. In the case
of classical knots (knots in 3-space), much of the theory concerns
such invariants derived from abelian coverings of the knot
complement. The algebraic information of these coverings
is contained in the Alexander module, the homology group of the
infinite cyclic cover, on which the deck transformation acts.
This module carries a bilinear form, the Blanchfield pairing.
Despite being previously well known, the interest of classical
abelian invariants is that only they remain easily and generally
computable, and therefore practically useful, despite their
various recently proposed (non-abelian \cite{Cochran,KL} and
quantum \cite{GK}) modifications. 

In this paper we will study some properties of algebraic
invariants for slice knots, i.e. knots trivial in the
topological concordance group. Our aim will be to find,
partly by computation, knots that demonstrate the failure
of possible implications between sliceness obstructions involving
various algebraic invariants. We believe that it is useful to
have concrete knots in hand to illustrate the occurring phenomena,
even though theory may suggest the existence of such examples.
 The recent expansion of knot tables \cite{HTW}
and computational tools has led to a series of such examples
related to other questions. In the present context, we are
mainly concerned by the undue lack of computations, even for
the (presumably easy to handle) abelian invariants.
In the final section, we will turn to smooth concordance, and
give some related examples.

\subsection{Topological concordance}

Classical (topological) knot concordance was introduced by Milnor and
Fox \cite{MilFox} in the 1950s. The decision whether two knots are
equivalent in that sense
is a longstanding problem, related to singularity theory, and
the classification of topological four-manifolds. Levine \cite{Levine}
made substantial early progress in the late 1960s, by introducing
an algebraic structure called the algebraic knot concordance
group, and using it to solve the problem in dimension at least 4.
While the algebraic concordance group is of fundamental importance
also in dimension 3, Casson and Gordon \cite{CG} showed that Levine's
homomorphism from the classical to the algebraic concordance group
has a non-trivial kernel.

The algebraic knot concordance group $C$ has a description within
a Witt group of quadratic forms (see for example \cite{Scharlau},
mainly \S 5). Levine's approach is to consider the Witt group
of isometric structures of the symmetrized Seifert form over $\bQ$.
(Working over $\bZ$ is much harder, and the arising invariants
were given later by Stoltzfus \cite{Stoltzfus}.) The Witt group
splits along primes $p$ in the ring $\bQ[t]$ of polynomials,
which are characteristic polynomials of the isometric structure.
Embedding $\bQ\subset \bR\subset \bC$, we obtain an integer invariant 
(signature) for polynomials $p$ (regarded now in $\bC[t]$)
with zeros on the complex unit circle $S^1$. Factoring out
$p$ (i.e. working over $\bQ[t]/\br{p}$), we have two $\bZ_2$ valued
invariants, the discriminant and Hasse-Witt invariant (latter
is not a homomorphism). Levine shows that the vanishing of all  
these invariants implies the vanishing of the form in the Witt group
(Theorem 21 in \cite{Levine3}) and that finite order elements in
$C$ have order 1, 2, or 4 (Proposition 22 (b) ibid.).

Now Levine defines, in \cite{Levine}, further algebraic knot
concordance invariants by Tristram-Levine signatures, and the
Alexander polynomial modulo Milnor-Fox factorizations. These must
thus find themselves in the above set of complete invariants of
$C$. In \cite{Ma} the relation is shown between Tristram-Levine
signatures and the real Witt group signatures. So by
Proposition 22 (a) in \cite{Levine3} invariants
detecting the (infinite rank) torsion-free part of $C$ correspond
exactly to Tristram-Levine signatures. The Milnor-Fox
condition relates to the discriminants (see proposition 5
in \cite{Levine}), and so if it holds, the algebraic concordance
class of the knot is completely determined by the Hasse-Witt
invariants. In theory these invariants can be realized
non-trivially, but no examples seem to have been elaborated on.

\subsection{Linking pairings and metabolizers\label{iqw}}

A different description of $C$ is in terms of Witt classes of
the Blanchfield pairing (valued in $\bQ[t]$), modulo pairings
with a self-annihilating submodule. There is a correspondence
between this and Levine's approach, first proven by Cherry
Kearton \cite{Kearton}. (Another source is the appendix A
of Litherland's note \cite{Litherland}, but beware of diverse
typos in the published version.) Under the
restriction $\bZ\to\bZ_n$ for any prime $n$, the Alexander module
turns into the homology group of the $n$-fold branched cyclic cover,
and the Blanchfield pairing determines a (non-singular) 
linking form on the torsion homology group of this cover.
In this paper we will be concerned with the case $n=2$.


Let $K$ be a knot, $D_K$ be its double branched cover, and $\lm$
the linking form on its $\bZ$-homology group $H_1=H_1(D_K)$
\cite{Lickorish}. The (finite and odd) order of $H_1(D_K)$ is
called the determinant $\dt=\dt(K)$ of $K$. The quadratic form
$\lm$ takes values in $\bZ/\dt(K)$, which is identified with the
subset of $\bQ/\bZ$ of fractions with denominator (dividing) $\dt(K)$.

If $K$ is slice (bounds a topological locally-flat disk in $B^4$),
then it is algebraically slice, i.e. all algebraic concordance
invariants of $K$ vanish. This occurs if and only if the Blanchfield
pairing has a self-annihilating submodule. So if the knot is
(algebraically) slice, then the $n$-fold cover linking forms are 
metabolic as well. Specifically for $n=2$ this means that then
$\lm$ vanishes on a subgroup $M$ of $H_1(D_K)$ of order
$\sqrt{\dt}$, equal to its annihilator.
$M$ is called a \em{metabolizer}. That $\dt(K)$ ought to be a square
is well-known from the condition of Milnor-Fox \cite{MilFox} that the
Alexander polynomial is of the form $\Dl_K(t)=f(t)f(1/t)$ for some
$f\in\bZ[t]$, since $\dt(K)=|\Dl_K(-1)|$. A further condition for $K$
being slice is that its (Murasugi) signature $\sg$ \cite{Murasugi}
vanishes, and so do the generalized (or Tristram-Levine)
signatures $\sg_\xi$, when $\xi$ is a unit norm complex number
and $\Dl(\xi)\ne 0$. (We have $\sg=\sg_{-1}$.)

The metabolizer existence condition is useful in some theoretical
situations, where the calculation of other invariants is more
tedious. See for example \cite{Livingston,Livingston2}. The
present work is mainly motivated by the interest in concrete
examples showing that this criterion is essential, in particular
as opposed to the other conditions for sliceness. We also investigate
the size of the isotropic cone $\Lm_0$ of the linking form.
We will find, in \S\reference{S4}, examples illustrating possible
phenomena concerning $\Lm_0$. First we give in \S\reference{S4.1}
computational examples, obtained from the tables of
\cite{KnotScape,HTW}, that show 

\begin{theorem}\label{theo1}
For each one of the three conditions below, there exist knots
satisfying this condition, which have zero Tristram-Levine
signatures and Bennequin numbers, and an Alexander polynomial
of the Milnor-Fox form.
\def\theenumi{(\alph{enumi})}
\def\labelenumi{\theenumi}
\begin{enumerate}
\item\label{(a)} $\Lm_0$ is trivial, i.e. $\{0\}$,
\item\label{(b)} $1<|\Lm_0|<\sqrt{\dt}$, 
\item\label{(c)} $|\Lm_0|\ge \sqrt{\dt}$,
  but $\Lm_0$ contains no subgroup of order $\sqrt{\dt}$.
\end{enumerate}
\end{theorem}

(Here `$\dt$' refers to the above introduced number $\dt(K)=|\Dl
_K(-1)|$. The meaning of Bennequin numbers is related to smooth
sliceness and explained in \S\reference{S5}.)


In \S\reference{S4.2} we find, now applying more systematical
constructions, examples that refine theorem \reference{theo1}.
We make decisive use of the realization of any admissible Alexander
polynomial by an unknotting number one knot. (This result was
proved first by Sakai \cite{Sakai} and Kondo \cite{Kondo}, and later
by several other authors, with a very recent construction due to
Nakamura \cite{Nakamura}.)

\begin{theorem}\label{theo2}
For any of the properties \reference{(a)}-\reference{(c)}
in theorem \reference{theo1}, we can find knots whose $H_1$ has
additionally no $4k+3$ torsion (or whose determinant
is not divisible by $4k+3$).
\end{theorem}

We can in fact classify all trivial cone forms on such groups $H_1$
(theorem \reference{t23}), although it is not clear which forms are
indeed realizable by knots.

The motivation for excluding $4k+3$ torsion lies in the structure of
the Witt group of $\bZ_p$-forms, and the possibility to rule out
concordance order two, when there is such torsion. For $4k+3$
torsion, recent further-going work of Livingston-Naik \cite{LN,LN2},
gives in fact sufficient conditions on infinite concordance order.
In \S\reference{S4.3} we will give examples where the isotropic
cone can exclude concordance order two, but Livingston-Naik's
criterion does not apply.

\begin{theorem}\label{theo3}
For the properties \reference{(b)} and \reference{(c)}
in theorem \reference{theo1}, we can find knots $K=\hK\# \hK$
so that the $p$-Sylow subgroup of $H_1(D_{\hK})$ for any prime
$p=4k+3$ is not cyclic of odd $p$-power order.
\end{theorem}

For property \reference{(a)} such knots do not exist. This follows
from an exact description of trivial cone forms on $\hK\#
\hK$ given in proposition \reference{t24} (which is similar to theorem
\reference{t23}, although all occurring forms are easily realizable).

\subsection{Smooth concordance invariants}

Knots are smoothly slice if they bound a smooth disk in
four-space. In the 1980s, Andrew Casson, using deep results of 
Freedman and Donaldson, gave the first example of a topologically
slice knot which is not smoothly slice. Such knots known by now
remain scarce (they all have trivial Alexander polynomial),
despite that some new candidates are recently suspected. The
difficulty in exhibiting such examples clearly displays the problems
with the methods we apply to study both types of concordance.

In recent years, the separation between topological and smooth
concordance grew wider with a vast development of new techniques
in the smooth category. One such is the inequality of Rudolph-%
Bennequin \cite{Bennequin,Rudolph2,Rudolph3}, which emerged in the
early 1990s, and sometimes proves (smooth) non-sliceness. More
recently, Ozsvath and Szabo \cite{OS} used Floer homology to define
an invariant $\tau$ that detects some non-slice knots. This invariant 
behaves similarly to the (Murasugi) signature and simultaneously
improves upon the Rudolph-Bennequin inequality. The Ozsvath-Szabo
invariant, in turn,
motivated Rasmussen \cite{Rasmussen} to define a (conjecturedly
equivalent) signature-like invariant $s$ from Khovanov homology.

In \S\reference{S5} we show how the inequality of Rudolph-Bennequin
can prove that a knot $K$ is non-slice, by applying this inequality
on knots $K'$ different from $K$. We also discuss the relation to the
recent knot homological ``signatures'' of Ozsvath-Szabo-Rasmussen, and
their status in the examples of \S\reference{S4}.

In \S\reference{SFi}, we explain concludingly how to construct prime
(in fact, hyperbolic and of arbitrarily large volume) knots
with the previously chosen properties.

\section{Preliminaries and notation\label{PS}}


\subsection{Knots, linking form and sliceness}

In the following knots and links will be assumed
oriented, but sometimes orientation will be irrelevant.

For a knot $K$, its \em{obverse}, or mirror image $!K$, is obtained
by reversing the orientation of the ambient space. The knot $K$
is called \em{achiral} (or synonymously \em{amphicheiral}),
if it coincides (up to isotopy) with its mirror image, and \em{chiral}
otherwise. When taking the knot orientation into account, we
write $-K$ for the knot $K$ with the reversed orientation. We
distinguish among achiral knots $K$ between $+$achiral and $-$achiral
ones, dependingly on whether $K$ is isotopic to $!K$ or $!-K$.

Prime knots are denoted according to \cite[appendix]{Rolfsen} for
up to 10 crossings and according to \cite{KnotScape} for $\ge 11$
crossings. We number non-alternating knots after alternating ones.
So for example $11_{216}=11_{a216}$ and $11_{484}=11_{n117}$.
We write $K_1\#K_2$ for the connected sum of $K_{1}$ and $K_{2}$,
and $\#^kK$ for the connected sum of $k$ copies of $K$.

Let $G$ be a finite group and let $p$ be a prime. A \em{Sylow
$p$-subgroup} of $G$ is a subgroup $H$ such that $p$ does not divide
$|G|/|H|$.
If $G$ is abelian, the \em{$p$-primary component} of $G$, written
$G_p$, is the subgroup of all elements whose order is a power of $p$.
It is the unique Sylow $p$-subgroup of $G$. By the \em{$p$-torsion
subgroup} we mean the subgroup of all elements whose order is equal
to $p$. 

By $D_K$ we denote the \em{double branched cover} of $S^3$ over a
knot $K$. (See \cite{CG,Rolfsen}.) By $H_1=H_1(D_K)=H_1(D_K,\bZ)$
we denote its homology group over $\bZ$. (The various abbreviated
versions will be used at places where no confusion arises; $H_1$
will be used throughout the paper only in this context, so that,
for example, when we talk of $H_1$ of a knot, always $H_1$ of its
double cover will be meant.) $H_1$ is a finite commutative group
of odd order. This order is called the \em{determinant} of a knot
$K$, and it will be denoted as $\dt=\dt(K)$. (It generalizes to
links $L$, by putting $\dt(L)=0$ to stand for infinite $H_1(D_L)$.)
By the classification of finite commutative groups, $H_1$ decomposes
into a direct sum of finite (odd order) cyclic groups $\bZ_k=\bZ/k=
\bZ/k\bZ$; their orders $k$ are called \em{torsion numbers}.

$H_1(K)$ is also equipped with a
bilinear form $\lm\,:\,H_1\times H_1 \to\bQ/\bZ$, called the
\em{linking form} (see \cite{Lickorish,MurYas} for example).
Since $\lm$ in fact takes values of the form $n/\dt$ for $n\in\bZ$,
we can identify them with $\bZ_{\dt}$.

A knot $K$ is called \em{slice} if it bounds a disk in $B^4$.
Except for \S\reference{S5}, and unless pointed out explicitly
otherwise, we work in the topological category. In \S\reference{S5},
we will consider smooth sliceness.

It is known that if $K$ is topologically slice, then $\lm$ is
\em{metabolic}. This means, there is a subgroup $M$ of $H_1(D_K)$
of order $\sqrt{\dt}$, which is equal to its annihilator 
\[
M^\perp\,=\,\{\,
g\in\,H_1\,:\,\lm(g,h)=0\,\mbox{ for all }\,h\in M\,\}\,.
\]
$M$ is called a \em{metabolizer}. Whenever
$\lm$ is non-degenerate, we have for any subgroup $G$ of $H$ that 
\begin{eqn}\label{Hp}
|G^\perp|\cdot |G|\,=\,|H|\,.
\end{eqn}

We also recall a few basic facts from number theory we will require
in the study of $\lm$. We apply for example Dirichlet's theorem on
infinitely many primes contained in arithmetic linear progressions.
We use also that every odd number $n$ is the sum of two squares, if
any only if every prime $p\equiv 3\bmod 4$ has even multiplicity $2e$
as factor of $n$, and then $p^{e}$ divides $a$ and $b$ in any solution
of $a^2+b^2=n$. Such facts can be found in standard books on number
theory; my personal favorites are \cite{HardyWright,Zagier}.

\subsection{Knot polynomials and signatures\label{Skp}}

The \em{skein polynomial} $P$ (introduced in \cite{HOMFLY}; here
used with the convention of \cite{LickMil}, but with $l$ and
$l^{-1}$ interchanged) is a Laurent polynomial in two
variables $l,m$ of oriented knots and links, and can be defined
by being $1$ on the unknot and the \em{(skein) relation}
\begin{eqn}\label{1} 
l^{-1}\,P\left(
L_+
\right)\,+\,
l \,P\left(
L_-
\right)\,=\,
-m\,P\left(
L_0
\right)\,.
\end{eqn}
Herein $L_{\pm,0}$ are three links with diagrams differing only near
a crossing. 
\begin{eqn}\label{Lpm0}
\begin{array}{*2{c@{\qquad}}c}
\diag{9mm}{1}{1}{
\picmultivecline{-8 1 -1.0 0}{1 0}{0 1}
\picmultivecline{-8 1 -1.0 0}{0 0}{1 1}
}
&
\diag{9mm}{1}{1}{
\picmultivecline{-8 1 -1 0}{0 0}{1 1}
\picmultivecline{-8 1 -1 0}{1 0}{0 1}
}
&
\diag{9mm}{1}{1}{
\piccirclevecarc{1.35 0.5}{0.7}{-230 -130}
\piccirclevecarc{-0.35 0.5}{0.7}{310 50}
}
\\[2mm]
L_+ & L_- & L_0
\end{array}
\end{eqn}
We call the crossings in the first two fragments respectively
\em{positive} and \em{negative}, and a crossing replaced by the third
fragment \em{smoothed out}. A triple of links that can be represented
as $L_{\pm,0}$ in \eqref{Lpm0} is called a \em{skein triple}. The
sum of the signs ($\pm 1$) of the crossings of a diagram $D$
is called the \em{writhe} of $D$ and written $w(D)$. The smoothing
of all crossings of $D$ yields the \em{Seifert circles} of $D$;
each crossing in $D$ can be viewed as connecting two Seifert circles.
Let $s(D)$ be the number of Seifert circles of $D$, and $s_-(D)$
be the number of those circles to which only negative crossings are
attached. We call such Seifert circles \em{negative Seifert circles}.

The substitution $\Delta(t) = P(-i,i(t^{1/2}-t^{-1/2}))$ (with $i=
\sqrt{-1}$) gives the (one variable) \em{Alexander
polynomial} $\Dl$, see \cite{LickMil}. It allows to express the
determinant of $K$, as $\dt(K)=\big|\,\Dl_K(-1)\,\big|$. The
(possibly negative) minimal and maximal power of $l$ occurring
in a monomial of $P(K)$ is denoted $\md_lP(K)$ and $\Md_lP(K)$.
Alexander polynomials (and factors thereof) will sometimes
be denoted by parenthesized list of their coefficients, putting
the absolute term in brackets. As an example of such a notation,
$1/t-1-t^2=(\ 1\ [-1]\ 0\ -1)$.

The \em{signature} $\sg$ is a $\bZ$-valued invariant of knots and
links. Originally it was defined in terms of Seifert matrices
\cite{Rolfsen}. We have that $\sg(L)$ has the opposite parity to the
number of components of a link $L$, whenever the determinant of $L$
is non-zero (i.e. $H_1(D_L)$ is finite). This in particular always
happens for $L$ being a knot, so that $\sg$ takes only even values
on knots.

Most of the early work on the signature was done by
Murasugi \cite{Murasugi}, who showed several properties
of this invariant. In particular the following property
is known: if $L_{\pm,0}$ form a skein triple, then
\begin{eqnarray}
\sg(L_+)-\sg(L_-) & \in & \{0,1,2\}\,, \label{2a} \\[2mm]
\sg(L_\pm)-\sg(L_0) & \in & \{-1,0,1\}\,. \label{2b}
\end{eqnarray}
(Note: In \eqref{2a} one can also have $\{0,-1,-2\}$ instead
of $\{0,1,2\}$, since other authors, like Murasugi, take $\sg$ to
be with opposite sign. Thus \eqref{2a} not only defines a property,
but also specifies our sign convention for $\sg$.) We remark that
for knots in \eqref{2a} only $0$ and $2$ can occur on the right. 

Let $M$ be a Seifert matrix for a knot $K$, and $\xi\in S^1$ a unit
norm complex number ($S^1$ denoting the set of such complex numbers).
The \em{Tristram-Levine} (or generalized) \em{signature} $\sg_\xi(K)$
of $K$ is defined as the
signature of the (Hermitian) form $M_\xi=(1-\xi)M+(1-\bar\xi)M^{T}$,
where bar denotes complex conjugation, and ${\,\cdot\,}^T$ means
transposition. We call $\xi$ and $\sg_\xi$ \em{non-singular} if their
corresponding form $M_\xi$ is so, that is, $\dt(M_\xi)\ne 0$, which
is equivalent to $\Dl_K(\xi)\ne 0$. For a fixed knot $K$ we obtain a
function $\sg_{*}(K)\,:\, S^1\to \bZ$ given by $\xi\mapsto \sg_\xi(K)$.
It is called the Tristram-Levine signature function of $K$. We
have $\sg=\sg_{-1}$, so Murasugi's signature is a special value
of $\sg_{*}$. If $\sg_\xi(K)$ is non-singular, then it is even.
(Since a knot has non-zero determinant, Murasugi's signature is
always non-singular.) Also $\sg_{*}(K)$ is locally constant around
non-singular $\xi$, that is, it changes values (``jumps'') only in
zeros $\xi$ of $\Dl_K$. The properties \eqref{2a} and \eqref{2b}
hold also for $\sg_\xi$.

Signatures (at least all
those we talk about in this paper) change sign under mirroring and
are invariant under orientation reversal, and so vanish on
amphicheiral knots. They are also additive under connected sum.

Let $g_t(K)$ be the topological \em{4-ball genus} of a knot $K$. Then
it is known, by Tristram-Murasugi's inequality, that if $\sg_\xi$ is
non-singular, then 
\begin{eqn}\label{TMi}
|\sg_\xi(K)|\le 2g_t(K)\,.
\end{eqn}
So if $K$ is topologically slice (that is, $g_t(K)=0$),
all non-singular $\sg_\xi$ vanish. Since a concordance
between $K_{1}$ and $K_{2}$ is equivalent to the
sliceness of $K_1\#-!K_2$, this implies that $\sg_{*}(K)$ is a
topological concordance invariant \em{outside} the zeros
of the Alexander polynomial. That is, if $K_{1,2}$ are concordant,
and $\Dl_{K_1}(\xi)\ne 0\ne \Dl_{K_2}(\xi)$, then $\sg_\xi(K_1)=
\sg_\xi(K_2)$. (In general one cannot say much about the behaviour
of singular $\sg_\xi$ under concordance \cite{Levine2}.)

Now more sophisticated methods are available to obstruct sliceness in
certain cases, like Casson-Gordon invariants \cite{CG} and twisted
Alexander polynomials \cite{Wada,KL}. Indeed, the Milnor-Fox and
Tristram-Murasugi conditions can be generalized to signatures and
twisted Alexander polynomials of certain non-abelian representations
of the knot group \cite{KL}. The general computability of such
invariants is still difficult, though (see \cite{Tamulis}). A similar
disadvantage reveal, for \em{smooth} sliceness, the very recent knot
homological (signature-like) concordance invariants $\tau$ and $s$ of
Ozsvath--Szabo and Rasmussen. The determination (or estimation)
of these invariants is often easier indirectly, using their properties,
rather than their definition. (We will later make a comment on
their calculation in relation to the Rudolph-Bennequin inequality.)

We invite the reader to consult \cite{unkn1,deg,4gen}
for more on the use of notation and (standard) definitions.

\section{The metabolizer criterion\label{S4}}

\subsection{Initial observations and remarks\label{S3.1}}

When looking for a metabolizer $M$ of $\lm$, there are some
simple, but important, observations to make.

First, $M$ always exists, whenever $H_1$ is cyclic. This already
restricts the search space for interesting examples, since about $80\%$
of the prime $\le 16$ crossing knots with square determinant have
cyclic (or trivial) $H_1$. It also suggests why composite knots
(where $H_1$ is more often non-cyclic) are likely to be of interest.

The second, and for us more relevant, observation is that
clearly any element in $g\in M$ must have $\lm(g,g)=0$. That is,
each $M$ is contained in the \em{isotropic cone}
\[
\Lm_0\,:=\,\{\,g\in H_1(D_K)\,:\,\lm(g,g)=0\in \bQ/\bZ\,\}\,
\]
of $\lm$. Note that if the isotropic cone contains a subgroup $G$,
then $G$ is always isotropic (the linking form is zero on $G$);
this is a bit more than a tautology, but follows from the identity 
\begin{eqn}\label{gh}
2\lm(g,h)\,=\,\lm(g+h,g+h)-\lm(g,g)-\lm(h,h)\,,
\end{eqn}
since we have no 2-torsion in $\bZ_{\dt}$. Thus a natural way to find
(or exclude the existence of) $M$ is to determine $\Lm_0$ and seek for
subgroups of $H_1$ of order $\sqrt{\dt}$ contained in $\Lm_0$. Any
such subgroup is a metabolizer (that is, equal to its annihilator)
because of \eqref{Hp}.

A third observation is that the linking form $\lm$, restricted to the
$p$-torsion subgroup of $H_1$ for a prime $p$, modulo metabolic
forms, naturally defines an element in the Witt group
of nonsingular $\bZ_p$ forms. This group is either $\bZ_2\op
\bZ_2$ or $\bZ_4$, depending on whether $p \equiv 1$ or $3 \bmod 4$ 
(see lemma 1.5, p.\ 87 of \cite{MH}). Thus, if $4k+3$-torsion
exists in $H_1(\hK)$, one may detect $\hK$ being 4-torsion
in the algebraic concordance group. The other invariants
(signatures, Alexander polynomial, etc.) can only detect (up to)
2-torsion. This suggests that if $\sg(\hK)=0$ and $4k+3\mid\dt
(\hK)$, then $K=\hKhK$ may be an example whose non-sliceness is
only detectable by the linking form. Therefore, in the search for
interesting examples $K=\hKhK$, we are led to consider (prime) knots
$\hK$ with $\sg=0$ and determinant divisible by $4k+3$.

In the case of $4k+3$-torsion, there is, though, Livingston
and Naik's result \cite{LN}, that if a prime $p=4k+3$ has
single multiplicity as divisor in $\dt$, then the knot has
infinite order in the classical knot concordance group. In
\cite{LN2} Livingston and Naik generalized their result to groups
$H_1$ whose Sylow $p$-subgroup is cyclic of odd power order.
Therefore, for interesting non-slice examples $K=\hKhK$, we
should consider in particular knots $\hK$ with determinant
divisible by a prime $p=4k+3$, but whose Sylow $p$-subgroup
of $H_1$ is not of the stated type for any such $p$.

For the computational part in the examples in \S\ref{S4},
we applied the computer program for calculating the
linking form \cite{unkn1}. This program was written in C
originally by Thistlethwaite, and later extended by myself. It
calculates the torsion numbers of $H_1$ and the corresponding
generators out of a Goeritz matrix \cite{GL} of a knot diagram.
The further algebraic processing was done with MATHEMATICA\TM{}
\cite{Wolfram}. Still computation alone often did not suffice
to find proper examples, and we were led to argue about their
\mbox{(non-)existence} mathematically.

{}From the explanation in \S\ref{iqw} is it clear that Witt group
invariants of the forms on all $n$-fold covers could tell something
about (algebraic) non-sliceness. But, again, we are unaware of any
concrete computations, even for the very restricted case $n=2$ we
consider. The study of higher $n$, from such computational point of
view, is certainly also worthwhile, and may be a future project.
(It is the lack of description in terms of Goeritz matrices of
higher cyclic cover homology groups that prevented from their
study.) It seems unclear (and may not be true) that the Witt class
of all finite cover forms recovers the Witt class of the
Blanchfield pairing.

\subsection{Examples with trivial, small or large cone\label{S4.1}}

The first series of examples shows that the existence of $M$ is
essential as opposed to the previously mentioned conditions on
Alexander polynomial, Rudolph-Bennequin numbers, and signature.

\begin{figure}[htb]
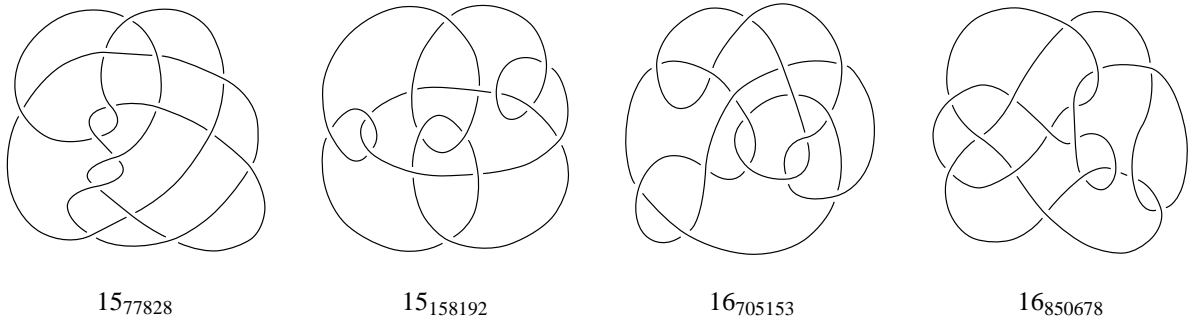

\[
\begin{array}{c*3{@{\qquad}c}}
\epsfsv{3.5cm}{15_77828} &
\epsfsv{3.35cm}{15_158192} &
\epsfsv{3.4cm}{16_705153} &
\epsfsv{3.4cm}{16_850678}
\\[19mm]
15_{77828} & 15_{158192} & 16_{705153} & 16_{850678} \\
\end{array}
\]
\caption{Knots with Milnor-Fox condition on the Alexander polynomial and
signature 0, but with trivial isotropic cone of the linking form.
The last one, $16_{850678}$, is $-$achiral (and fibered of genus 4).
\label{fig1}}
\end{figure}

\begin{exam}\label{exam1}
Consider the knots on figure \reference{fig1}. The first (alternating)
knot $15_{77828}$ has $\sg=0$ and Alexander polynomial 
\[
\Dl=(\ 1\ -8\ 32\ -82\ 152\ -216\ [243]\ -216\ 152\ -82\ 32\ -8\ 1)\,,
\]
which is of the Milnor-Fox form $f(t)f(1/t)$ with
$f=(\ [1]\ -4\ 8\ -9\ 8\ -4\ 1)$. (We explained in the appendix of
\cite{deg} that for given $\Dl$ only finitely many $f$ come in
question, and one can make the search for $f$ very efficient. Thus
the Milnor-Fox test is easy to perform.) The Rudolph-Bennequin
inequality (in the \em{smooth} setting; see \S\reference{S5}, and
in particular remark \ref{r5.1}) is also trivial on its 15 crossing
diagrams. This knot has $H_1=\bZ_{35}\op\bZ_{35}$. Two particular
generators $g_{1}=(1,0)$ and $g_{2}=(0,1)$ of the cyclic factors have
\[
\lm(g_1,g_1)=8/35,\ \ 
\lm(g_1,g_2)=18/35,\ \ 
\lm(g_2,g_2)=4/35.\ \ 
\]
This shows that $\Lm_0=\{(0,0)\}$. Hence $15_{77828}$ is not slice.
The knots $15_{158192}$, $16_{705153}$, $16_{747143}$ and $16_{850678}$,
three of which also appear in figure \reference{fig1}, are of similar
nature. They all have $\Dl=(\ 1\ -10\ 43\ -100\ [133]\ -100\ 43\ -10\ 
1)$ of Milnor-Fox form, $H_1=\bZ_{21}\op\bZ_{21}$, and trivial $\Lm_0$.
Particularly interesting is $16_{850678}$, because it is $-$achiral.
Thus for it all (including the singular) Tristram-Levine signatures
vanish. (Similarly, as will follow from the explanation in \S\ref{S5},
the Rudolph-Bennequin method fails~-- with certainty; quite likely
it also fails for the others.)
\end{exam}


\begin{figure}[htb]
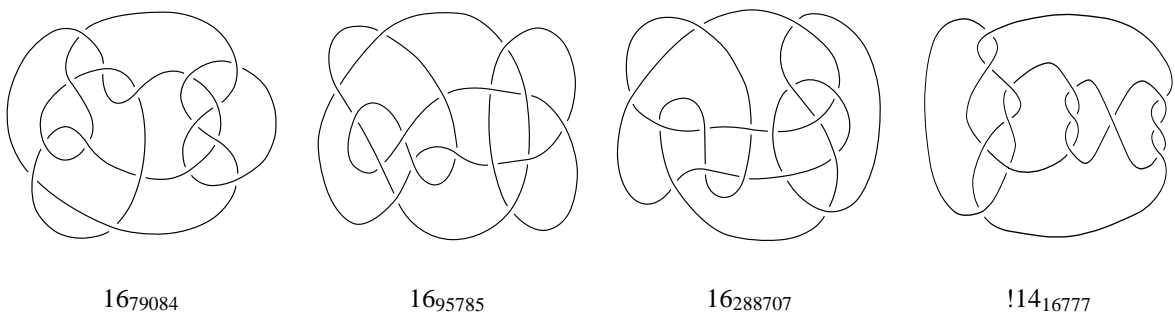

\[
\begin{array}{c*3{@{\enspace\quad}c}}
\epsfsv{3.3cm}{16_79084} &
\epsfsv{3.15cm}{16_95785} &
\epsfsv{3.2cm}{16_288707} &
\diag{1cm}{3.3}{3.0}{
\picscale{1 1.38}{
\picputtext[dl]{0 0}{%
\autoepsffalse\epsfs{3.3cm}{k-14_16777}%
}}}
\\[19mm]
16_{79084} & 16_{95785} & 16_{288707} & !14_{16777} \\
\end{array}
\]
\setbox\@tempboxa=\hbox{$\sqrt{\dt}$}
\caption{Knots providing examples with 
the Milnor-Fox condition on the Alexander polynomial and
signature 0, and with non-trivial isotropic cone of the linking form,
which still does not contain a metabolizer. The first three
have $|\Lm_0|=3$, while the last one's connected sum with $3_1\# 3_1$
has $|\Lm_0|\ge \box\@tempboxa$.\label{fig2}}
\end{figure}

\begin{exam}
The fact that in the above series $\Lm_0$ is always trivial originally
led to suspect that when Milnor-Fox holds and $\sg=0$, then $\Lm_0\ne 0$
may already imply $\Lm_0\supset M$. There is no reason why this should
be true, but the examples found to refute it required a considerable
quest. 
There were a total
of 10 prime knots of $\le 16$ crossings (all alternating 16 crossing
knots),
for which determinant is a square, $\sg=0$ and $\Lm_0\ne 0$,
but $\Lm_0\not\supset M$. Seven of them (three of which are shown
in figure \reference{fig2}, and are the knots to the left) have
$\Dl(t)=f(t)f(1/t)$. For all 10 knots $H_1=\bZ_{135}\op\bZ_{15}$,
with $\Lm_0=\{(0,\ 0),\ (45,\ 0),\ (90,\ 0)\}$, which is still
too small.
\end{exam}

Thus a further-going and more complicated question is what occurs if
we do not desire that $\Lm_0\supset M$ is excluded already because
of cardinality reasons, that is, if $|\Lm_0|\ge\sqrt{\dt}$. 

\begin{exam}\label{xyz}
Consider the knot $K=3_1\#3_1\#!14_{16777}$. (The knot $!14_{16777}$
is given on the right in figure \reference{fig2}.) We have
$\sg(14_{16777})=4$ and $\Dl(14_{16777})=\Dl(3_1)^2f(t)f(1/t)$,
where $f(t)=-2+2t-2t^2+t^3$. Clearly $\sg(K)=0$, but since $f(t)$
has no zeros on the unit circle, we can conclude that even all
non-singular (i.e. not corresponding to the zeros of $\Dl(3_1)$)
Tristram-Levine signatures of $K$ vanish. We have $\dt=3969$, and
$H_1=\bZ_{147}\op \bZ_3\op \bZ_3\op\bZ_3$. It turns out that
$|\Lm_0|=231>63=\sqrt{\dt}$, but $\Lm_0$ contains no subgroup
of order $63$.

A more complicated knot with simpler factors is $K'=5_2 \# 8_{18}
\# !9_{12}$. There we have
\[
\Dl(t)\,=\,(1-t+t^2)^2\cdot (1-3t+t^2)^2\cdot (2-3t+2t^2)^2/t^6\,.
\]
Since both $1-t+t^2$ factors come from (the Alexander polynomial
of) $8_{18}$, which is amphicheiral,
one can similarly conclude that all non-singular Tristram-Levine
signatures of $K'$ vanish. Now $H_1=\bZ_{105}\oplus \bZ_{105}$, and
$|\Lm_0|=117$, but $\Lm_0$ contains no subgroup of order $105$.
\end{exam}

\begin{remark}
The most general condition under which one can seek a large metabolizer,
of course, is when $\dt(K)$ is a square. (Otherwise, the definition
of a metabolizer does not make sense.) It turned out that even in this
most general setting, among \em{prime} knots up to 16 crossings, there
was a single example, and it does not fall into any of the further
specified categories. This example is $15_{197573}$. (This knot has
$\sg=-4$ and $\Dl(t)\ne f(t)f(1/t)$ and 15 crossing diagrams with
non-trivial Rudolph-Bennequin numbers.) It has $H_1=\bZ_3\op\bZ_3\op
\bZ_3\op\bZ_3$, and $|\Lm_0|=21$, but $\Lm_0$ contains no subgroup
of order $9$. 
\end{remark}

%

It is also interesting to ask what values $|\Lm_0|$ can attain. In
particular how large can $|\Lm_0|$ become for given determinant?
Can $|\Lm_0|$ be relatively prime to $\dt$?
Obviously it does not need to share all determinant's prime divisors.


\subsection{Examples with no $4k+3$-torsion\label{S31}\label{S4.2}}

For some explicitly named knots of $11$ or more crossings that occur
in the examples of \S\ref{S31} and \S\ref{S32}, see figure \ref{fig4}.

\begin{exam}\label{ex4}
The search among connected sums allows also to find knots $K$, where
$H_1$ has no $4k+3$-torsion. They were motivated by the remarks about the
Witt group in \S\ref{S3.1}. One least crossing number example we found
is $K=6_3 \# 12_{1152}$. (Here $H_1=\bZ_{65}\oplus \bZ_{65}$.)
One verifies similarly all Alexander polynomial and signature
conditions, but $|\Lm_0|=1$. 
\end{exam}

All other examples we found have trivial $\Lm_0$ as well.
However, there is a systematic way of constructing knots
with large $\Lm_0$, which we explain.

\begin{prop}
There are knots $K$ with $4k+3\nmid \dt$, with $\Dl=f(t)f(1/t)$
and zero Tristram-Levine signatures, with $|\Lm_0|>\sqrt{\dt}$
but no metabolizer.
\end{prop}

\proof
Let $d\equiv 1\bmod 4$ have no $4k'+3$ divisor, and let it contain
any prime with multiplicity one. Let $\Dl(t)\in\bZ[t,t^{-1}]$
be some polynomial with $\Dl(t)=\Dl(1/t)$, $\Dl(1)=1$, no zero on
the unit circle and $\Dl(-1)=d$. That such $\Dl$ exists is easy.
Consider the base-4-expansion 
\[
d\,=\,\sum_{i=0}^n\,e_i4^i
\]
of $d$ (with $0\le e_i\le 3$ and $e_0=1$), and take $\Dl(t)=\nabla (t^
{1/2}-t^{-1/2})$ where $\nabla(t)\,=\,\sum e_i(-t^2)^i$. By
Kondo's result \cite{Kondo}
there is an unknotting number one knot $K$ with $\Dl(K)=\Dl$. By
\cite{Wendt}, $H_1(K)$ is cyclic, and by \cite{Lickorish,unkn1}, there
is a generator $g$ of $H_1$ with $\lm(g,g)=2/d$ (note that $\sg(K)=0$). 

The prime condition assures that $H_1$ has no elements of
non-trivial prime power order, and then the metabolicity of $\lm$ is
equivalent to the metabolicity of its reductions on the $p$-torsion
subgroups, for all primes $p$ dividing $d$. Then by the Witt
group argument, $\lm\op\lm$ on $H_1\op H_1$ is metabolic. 

Since $\lm\op\lm$ is metabolic, $|\Lm_0(K^{\#2})|\ge \dt(K)=d$. Now,
because of \eqref{gh},
\[
\Lm_0(K^{\#4})\,\supset\,\Lm_0(K^{\#2})\op \Lm_0(K^{\#2})\,,
\]
and $|\Lm_0(K^{\#4})|\ge d^2$. We want to show now that this inequality
is strict, and so we must show that the inclusion is proper. Now, when
$d=4k+1$ has no $4k'+3$ divisors, there are $a$ and $b$ relatively
prime to $d$ with $a^2+b^2=d$. Then $(a,0,b,0)$ lies in
$\Lm_0(K^{\#4})$ but not in $\Lm_0(K^{\#2})\op \Lm_0(K^{\#2})$.

Then for any knot $K_0$ with non-metabolic $\lm$, the knot
$K_0\# K^{\#2k}$ will also have non-metabolic $\lm$, but for
$k$ large enough $|\Lm_0|>\sqrt{\dt}$. \qed

In the case of trivial cone, we can in fact classify the
forms algebraically in a slightly more general situation.
This explains the nature of the knots found computationally
in example \reference{ex4}.

\begin{theorem}\label{t23}
Assume $H$ is a finite commutative group of odd order $d$, and that
$d$ has no $4k+3$ divisors. Let $\lm\,:\,H\,\times\,H\,\to\,\bQ/\bZ$
is a symmetric bilinear form with trivial cone $\Lm_0$.
Then (and only then) $H=\bZ_{q'}\op \bZ_{q}$, where $q'$ is
a product of distinct primes $p_i=4l_i+1$, and $q\mid q'$, and for
each prime $p_i\mid q$, there is a basis $(g_1,g_2)$ of the
$p_i$-Sylow subgroup $\bZ_{p_i}\op\bZ_{p_i}$ of $H$ such that with
$(e_1,e_2)$ denoting $e_1g_1+e_2g_2$ we have
\begin{eqn}\label{st**}
\lm\bigl(\,(e_1,e_2),(e_1,e_2)\,\bigr)\,=\,e_1^2\,+\,\bt e_2^2\,,
\end{eqn}
where $\bt$ is not a square residue mod $p_i$.
\end{theorem}

\proof Of course if $\Lm_0=0$, then $\lm$ is non-degenerate.
We have $\Lm_0$ if and only if $[\Lm_0]_{p_i}=0$, where
$[\Lm_0]_{p_i}$ is the reduction of $\Lm_0$ on the $p_i$-Sylow
subgroup (or $p_i$-primary component) $H_{p_i}$ of $H$. 

Fix a prime $p=p_i\mid d$. If $\bZ_{p^2}\subset H_{p}$,
then one easily finds non-trivial elements in $[\Lm_0]_{p}$.
Thus $H_{p}=\bZ_{p}^{\op k}$. If $k=1$ then $[\Lm_0]_{p}=0$,
since $\lm$ does not degenerate on $H_{p}$.
Consider then $k=2$. Then $H_{p}=\bZ_{p}\op\bZ_{p}$ with
a basis $(g_1,g_2)$. We write $(e_1,e_2)$ for $e_1g_1+e_2g_2$.
Then for some $a,b,c\in\bZ_p$ we have
\begin{eqn}\label{st*}
\lm\bigl(\,(e_1,e_2),(e_1,e_2)\,\bigr)\,=\,ae_1^2\,+\,be_1e_2\,+\,c
e_2^2\,.
\end{eqn}
We can assume $a\ne 0\ne c$, else $\Lm_0\ne 0$.
Then since $2$ is invertible in $\bZ_{p}$, we have
\[
\lm\bigl(\,(e_1,e_2),(e_1,e_2)\,\bigr)\,=\,
a\,\left(\,e_1+\frac{b}{2a}e_2\,\right)^2+\,
\left(c-\frac{b^2}{4a}\,\right)e_2^2\,.
\]
Now $\bigl(e_1\mapsto e_1+\frac{b}{2a}e_2,\,e_2\mapsto e_2\,\bigr)$
is bijective, and so we can assume w.l.o.g. that in \eqref{st*}
we have $b=0$. 

Now the multiplicative group $\bZ_p^*$ of units of $\bZ_p$ is cyclic.
Thus there are two equivalence classes of non-trivial residue classes
modulo $p$ up to multiplication with squares. If $a$ and $c$
lie in the same class then we can assume w.l.o.g. that $a=c$,
find $4k+1=p=e_1^2+e_2^2$, and have $\Lm_0\ne 0$. If $a$
and $c$ are in different classes, then one can make (exactly)
one of them equal $1$, and $\lm$ has the form in \eqref{st**}.

Now let $k\ge 3$. Consider $\bZ_{p}^{\op 3}\subset H_p$ with
a basis $(g_1,g_2,g_3)$. Then one can again assume that
$m_i=\lm(g_i,g_i)\ne 0$, and by substitutions diagonalize
$\lm$. That is, we can assume w.l.o.g. that $\lm(g_i,g_j)=0$
when $i\ne j$. Then, since at least two of $m_1,m_2,m_3$
are in the same equivalence class modulo squares in $\bZ_p^*$,
by the previous argument we can find a non-trivial element
in the cone. \qed

\begin{remark}\label{r2}
If now $\dt$ is a square, then for trivial $\Lm_0$ we must have $q'=q$.
Among the knots obtained from our calculations we had $q$ being the
product of two primes, $p_1=5$ and $p_2\in\{13,17,29,37\}$. (For each
of these $p_2$ examples are $4_1\#6_3\#12_{1722}$, $8_3\#14_{2160}$,
$11_{76}\# 11_{160}$ and $8_{17}\# 14_{15965}$ resp.) No knots occurred
where $q$ is only a single prime. More particularly, from 122,624
prime $\le 16$ crossing knots with determinant $d=4k+1$ prime and
$\sg=0$, all have $\lm(g,g)= 2/d$ for some generator $g$ of $H_1$.
(For any $d$ a knot $K_1$ of such linking form always exists, so
if another knot $K_2$ fails to realize this form, then $K_1\#K_2$ is
a potential candidate for trivial $\Lm_0$.) This phenomenon seems
related to some property of Minkowski units, but so far I cannot
work out an exact explanation.
\end{remark}

Now we explain how to construct knots with no $4k+3$-torsion
and $1<|\Lm_0|<\sqrt{\dt}$.

Let $K_2$ be a knot with determinant $d$ being a product of distinct
primes $p=8k+5$, with a generator of (the necessarily cyclic)
$H_1$ having $\lm(e_2,e_2)=1/d$. Let $K_1$ be a knot with
$\dt=d^3$ and cyclic $H_1$ with $\lm(e_1,e_1)=2/d^3$.

Then, for each prime $p\mid d$, consider $\Lm_0(K_1\# K_2)$
on the $p$-Sylow subgroup $\bZ_{p^3}\op\bZ_p$ of $H_1(K_1\# K_2)$.
We calculate $\big|[\Lm_0]_p\big|$. We have in $\bQ/\bZ$
\[
\lm\bigl(\,(e_1,e_2),(e_1,e_2)\,\bigr)\,=\,2e_1^2/p^3\,+\,e_2^2/p\,.
\]
If $p\nmid e_1$, then the first term on the right has denominator
$d^3$, and so $\lm\ne 0$. Thus $p\mid e_1$. Let $e_1'\in\bZ_p$ be
the reduction of $e_1/p\in \bZ_{p^2}$ modulo $p$. Then $e_1',e_2
\in\bZ_p$ satisfy $2e_1'^2+e_2^2=0$. Since $p\equiv 5\bmod 8$,
we have that $\pm 2$ is not a quadratic residue, and $e_1'=e_2=0$.
Thus the vectors $v=(e_1,e_2)$ with $\lm(v,v)=0$ are multiples
of $(p^2,0)$, and so $\big|[\Lm_0]_p\big|=p$.

Thus $|\Lm_0|=d$, while $|H_1|=d^4$. With this idea in mind
we can find examples.

\begin{exam}
Consider the knot $K_2=11_{333}$, which is also the 2-bridge
knot $(65,14)$. It has signature $\sg=0$, determinant $d=65$
and a generator $g$ of $H_1$ with $\lm(g,g)=14/65$, which is
equivalent to $1/d$ up to squares. (By remark \reference{r2},
this is apparently the smallest $d$ for which we can find
$K_2$.) The Alexander polynomial
\[
\Dl\,=\,(\ 4\  -16\   25\  -16\    4)
\]
has no zero on the unit circle. Let $K_1$ be a knot of unknotting
number one, whose Alexander polynomial $\Dl_{K_1}$
has no zero on the unit circle,
is of the form $\Dl_{K_2}f(t)f(1/t)$ and $\Dl_{K_1}(-1)=65^3=274,625$.
(For example take $\Dl_{K_1}=\Dl_{K_2}^3$.) Then $K_1\# K_2$ is a
knot of the type we sought.
\end{exam}

\subsection{Excluding concordance order 2\label{S32}\label{S4.3}}

Here we present some examples where our method prohibits the sliceness
of the connected sum of a knot with itself, i.e. rules out (topological)
concordance order 2. Among $\hKhK$ type examples, we can find or
construct the following knots.

\begin{exam}
The simplest example arising is $7_7\# 7_7$
(with $H_1=\bZ_{21}\oplus \bZ_{21}$). It has trivial $\Lm_0$.
This can be easily explained, since $7_7$ has unknotting number one,
and $a^2+b^2=0$ has no non-trivial solutions in $\bZ_{21}$.
\end{exam}

There are, however, no examples of trivial cone, to which the
Livingston-Naik result does not apply. This is explained followingly.

\begin{prop}\label{t24}
$\lm(\hK\# \hK)$ has trivial isotropic cone, if and only if $\dt(\hK)=d
$ is a product of distinct primes $p$, all of which congruent $3\bmod
4$.
\end{prop}

\proof Note that $\Lm_0(\hK\# \hK)$ is trivial, if and only if it is
so in every reduction to prime torsion subgroups of $H_1$.

If the determinant $d$ is of the described exceptional type,
then the reduction of $\Lm_0$ to $\bZ_p\op\bZ_p$ is trivial, since
$p$ is not of the form $a^2+b^2$. So $\Lm_0$ is trivial.

Assume now the determinant $d$ is not of the specified type. If a
prime $p=4k+1$ divides the determinant $d$, then the reduction
of $\Lm_0(\hK\# \hK)$ on the $p$-torsion subgroup is metabolic
by the Witt group argument, and so $\Lm_0$ cannot be trivial.

Now let $p=4k+3$ be a prime with $p^2$ dividing $d$.
If $\bZ_{p^2}$ occurs as subgroup, one immediately
finds non-trivial zero-linking elements. Thus the
$p$-Sylow subgroup of $H_1(\hK)$ must be a multiple of
$\bZ_p$, and since we assume $p^2\mid\dt$, we have at
least two copies of $\bZ_p$. Since $\lm(\hK)$ is
non-degenerate, each generator of each $\bZ_p$ has
non-zero linking. Now the multiplicative group $\bZ_p^*$
of units of $\bZ_p$ is cyclic, and so any element
is a plus or minus a square. Then on $\bZ_p\op
\bZ_p\subset H_1(\hK)$ we have in $\bZ_p$ up to sign
\begin{eqn}\label{pq}
\lm\left((e_1,e_2),(e_1,e_2)\right)=e_1^2\pm e_2^2+qe_1e_2
\end{eqn}
for some $q\in\bZ_p$.
Since $2$ and $4$ are invertible in
$\bZ_p$, we can write in $\bZ_p^{\op 4}\subset H_1(\hK\# \hK)$
\begin{eqn}\label{pq2}
\lm\bigl((x,y,z,w),(x,y,z,w)\bigr)\,=\,\left(x+
\frac{q}{2}y\right)^2+\left(z+\frac{q}{2}w\right)^2+q'(y^2+w^2)\,,
\end{eqn}
with $q'=\left.\pm 1-\myfrac{q^2}{4}\right.$.

If we have negative sign in \eqref{pq}, then $(x,x,x,-x)$ for every
$x$ is isotropic. So assume we have positive sign. If $q=\pm 2$, then
in \eqref{pq} we have $\lm=(e_1\pm e_2)^2$, and we are easily done.

So assume $q\ne \pm 2$. Then
$-q'$ is invertible in $\bZ_p$. Consider the arithmetic progression
$-1/q'+k'p$, and the subprogression in it made of numbers $4k+1$.
If $-1/q'$ is even, then $k'$ is odd and vice versa. Thus
we have a progression $a'+b'\cdot(2p)$, where $(a',2p)=1$. This
progression contains a prime $r$ by Dirichlet's theorem, and since
$r=4k+1$, we have $-1/q'+k'p=r=y^2+w^2$ for some $y,w$ (obviously
not both divisible by $p$, since $(p,r)=1$). Then in $\bZ_p$
we have $q'(y^2+w^2)=-1$, and for these $y$ and $w$ we can find
$x$ and $z$ with $x+\myfrac{q}{2}\,y=0$ and $z+\myfrac{q}{2}w=1$,
so by \eqref{pq2} we are done.
\qed


\begin{figure}[htb]
\[
\begin{array}{c*3{@{\qquad}c}}
\epsfs{2.95cm}{t1-11_274} & \epsfs{3.0cm}{t1-11_333_2}
& \epsfs{3.2cm}{t1-12_554_7} & \epsfs{3.1cm}{t1-12_1152} 
\\[4mm]
11_{274} & 11_{333} & 12_{554} & 12_{1152} \\
\end{array}
\]
\caption{\label{fig4}}
\end{figure}



\begin{exam}
Among small cone type examples, we found the knot $11_{274}\# 11_{274}$.
The Alexander polynomial of $\hK=11_{274}$ is 
\[
(\ 1\   -6\   18\  -35\   [45]\  -35\   18\   -6\    1)
\]
with determinant 165, and no zeros on the unit circle.
In this case $\hKhK$ has $|\Lm_0|=9$. A similar example is
$\hK=11_{280}$.
\end{exam}

To find examples, to which the Livingston-Naik results do not apply,
let $\hK=\hK_1\#\hK_2$, where $\hK_{1,2}$ are unknotting number one
knots, whose Alexander polynomial has no zero on the unit circle.
We choose the determinants of $\hK_{1,2}$ to be $d^2$ and $d$ resp.,
and $d$ to be a product of an even number of different primes
$p_i=4k_i+3$.

To calculate $|\Lm_0(\hK\# \hK)|$ it suffices to consider its
restriction on the $p$-Sylow subgroup for $p$ being any of the
$p_i$. This subgroup is $\bZ_{p^2}\op\bZ_{p^2}\op\bZ_{p}\op\bZ_{p}$.
Now if $p$ divides $a^2+b^2$, then so does $p^2$, and $p$ divides
both $a$ and $b$. This means first that the $p$-adic valuation of
the restriction of $\lm$ on $\bZ_{p^2}\op\bZ_{p^2}$ is 0 or $-2$,
so that $[\Lm_0]_p$ splits into a direct sum over its part in
$\bZ_{p^2}\op\bZ_{p^2}$ and $\bZ_{p}\op\bZ_{p}$. It means second
that former summand of $[\Lm_0]_p$ has size $p^2$, while latter
summand is trivial. Thus $\big|[\Lm_0]_p\big|=p^2$, and
$|\Lm_0|=d^2$.



\begin{exam}
We found among low crossing knots one single example of the large
cone type. Here $\hK=12_{554}$. This knot has $H_1=\bZ_{21}\op\bZ_3
\op\bZ_3$. The Alexander polynomial 
\[
\Dl=\ (\ -2\   15\  -45\   [65]\  -45\   15\   -2)
\]
has determinant 189 and no zeros on the unit circle. We have for
$K=\hKhK$ that $|\Lm_0|=225$, but $\Lm_0$ contains no subgroup
of order $189$.
\end{exam}

For an example, to which the Livingston-Naik result does not apply,
take $\hat K=7_7\# 12_{554}^{\# 2k}$. Since $H_1(12_{554})$ has no
elements of non-trivial prime power order, the Witt group argument
ensures that $\lm(\hKhK)$ is not metabolic, because it is not so
for $k=0$. But since $|\Lm_0(12_{554}^{\# 2})|>\dt(12_{554})$, we
obtain for large $k$ again large cone. 

\section{Indirect Rudolph-Bennequin inequality\label{S5}}

For the final contribution of the paper, we turn to \em{smooth}
sliceness. In contrast to the topological case in \eqref{TMi},
let $g_s(K)$ be the smooth $4$-genus of $K$. So $K$ is smoothly
slice if and only if $g_s(K)=0$.

Rudolph \cite{Rudolph2} showed the ``extended slice Bennequin
inequality'' (later proved also, and slightly more clarified
by Kawamura \cite{Kawamura}), which gives a lower estimate
for $g_s(K)$ from a reduced diagram $D$ of $K$:
\begin{eqn}\label{sbi}
g_s(K)\ge \frac{w(D)-s(D)+1}{2}+s_-(D)\,=:\,rb(D)\,.
\end{eqn}\unskip
As we explained in \S\ref{Skp}, by $w(D)$ we denote the writhe of
$D$, and by $s(D)$ and $s_-(D)$ the number of its Seifert circles resp.
negative Seifert circles. (Kawamura remarked that we must exclude
diagrams with negative Seifert circles adjacent to nugatory crossings.)
We call $rb(D)$ the \em{Rudolph-Bennequin number} of $D$. Inequality
\eqref{sbi} is an improvement of Bennequin's original
inequality \cite[theorem 3]{Bennequin}, which estimates the ordinary
genus $g(K)$ of $K$ by the \em{Bennequin number} 
\[
b(D)=\frac{w(D)-s(D)+1}{2}\,,
\]
in which the $s_-(D)$ term is missing. Rudolph showed prior to the above
improvement \eqref{sbi}, that $b(D)$ also estimates $g_s(K)$ (``slice
Bennequin inequality''). 

This quantity $b(D)$ has another upper bound,
namely the minimal degree $\md_lP(K)$ of the
skein polynomial $P$, as proved by Morton \cite{Morton}.
In particular, if 
\[
\dl(P(K))\,=\,\max\bigl(\md_lP(K),- \Md_lP(K)\bigr)\le 0\,,
\]
the original
Bennequin number $b(D)$ will be useless in showing that $K$ is not
slice, whatever diagram $D$ of $K$ (or its mirror image) we apply
it to. The $s_-(D)$ term, however, lifts the skein polynomial
obstruction, and in \cite{4gen} we showed that indeed $\md_lP=0$
can still allow the existence of diagrams with $rb>0$. (Our example
$13_{6374}$ has also $\Dl=1$, so that any other previous method to
prohibit sliceness fails.) 

In theory thus, for many non-slice
$K$, we could have $rb(D)>0$ for some $D$. In practice, however,
to find such $D$ for most $K$ is a tedious, or even pointless,
undertaking. It may well be that $D$ does not exist, and it is
not worth checking more than a few diagrams that can be easily
obtained. Since the improvement involving $s_-(D)$ is modest,
$\dl(P)\ll 0$, even if not a definite obstruction, still remains at
least a good heuristic evidence that such $D$ is unlikely to exist.

There is one particular situation, in which one can definitely
exclude the existence of $D$. Namely, note that it would imply
that
\begin{eqn}\label{iy}
\lim_{n\to\infty}\,g_s(\#^n\pm K)\,=\,\infty\,,
\end{eqn}
as $rb$ is additive under connected sum of diagrams (if it is properly
performed, and unless the estimate is trivial), and invariant under
reversal of knot orientation. Thus in particular if $K$ is of finite
order in (smooth) concordance \em{as an unoriented knot}, $D$ cannot
exist.
The finite concordance order amounts, at least in practice, to saying
that $K$ is slice or achiral (of either sign), explaining the special
role of $16_{850678}$ in example \reference{exam1}.

\begin{figure}[htb]
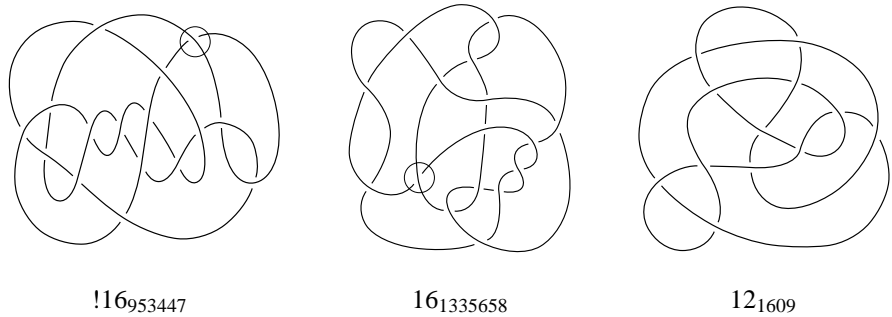

\[
\begin{array}{c*2{@{\qquad}c}}
\vcbox{
  \diag{1cm}{3.6}{3.35}{
    \picputtext[dl]{0 0.1}{\epsfs{3.35cm}{16_953447}}
    \piclinewidth{90}
    \piccircle{2.54 2.81}{0.2}{}
  }
} &
\vcbox{
  \diag{1cm}{3}{3.35}{
    \picputtext[dl]{0 -0.1}{\epsfs{3.35cm}{16_1335658}}
    \piclinewidth{90}
    \piccircle{0.96 1.0}{0.2}{}
  }
} &
\epsfsv{3.4cm}{12_1609} 
\\[19mm]
!16_{953447} & 16_{1335658} & 12_{1609} \\
\end{array}
\]
\caption{The two 16 crossing knots satisfy the Milnor-Fox condition on the
Alexander polynomial (the second one has trivial polynomial), have signature
0, and cyclic $H_1(D_K)$. On any of their 16 crossing diagrams that
could be found, the Rudolph-Bennequin inequality is also trivial.
However, in the diagrams depicted, a crossing switch results in a
diagram of $12_{1609}$. This knot has the diagram on the right, with
(Rudolph-)Bennequin number $2$. Thus the $4$-genus is at least 2, and
so the original 16 crossing knots are not slice.\label{fig3}}
\end{figure}

As a contrast to the intuition described so far, we
conclude by showing how the Rudolph-Bennequin inequality
can prove indirectly that some knots are
not slice, when the quest for a diagram $D$ with $rb(D)>0$
(among the suggestive candidates) fails. See figure \reference{fig3}.

\begin{exam}
Consider the knot $!16_{953447}$. Using Thistlethwaite's tools,
one can generate 253 different 16 crossing diagrams of this knot,
all of which, however, have $rb\le 0$. This knot has the
Alexander polynomial of the square knot, and $\sg=0$.
(This again implies that all Tristram-Levine signatures vanish.)
Its $H_1$ is cyclic, and hence $\lm$ is metabolic.

Now we pursue the following idea. If $K'$ differs from $K$ by one
crossing change (performed in whatever diagram $D$ of $K$), then
$|g_s(K')-g_s(K)|\le 1$, so that if $K$ is slice, $g_s(K')\le 1$.
However, if we find a diagram $D'$ of $K'$ with $rb(D')\ge 2$, then
$g_s(K')\ge 2$, and so $K$ can not be slice. Consider the diagram of
$K=!16_{953447}$ on the left of figure \reference{fig3}. Switching
the encircled crossing turns $K$ into $K'=12_{1609}$. When $K'$
is depicted in the diagram $D'$ on the right of figure \ref{fig3},
then $b(D')=rb(D')=2$, and so we can conclude indirectly that
$!16_{953447}$ is not slice. 
\end{exam}

\begin{exam}
The same argument (again with $K'=12_{1609}$) applies to $16_{1335658}$.
One can also handle in a similar way (using another knot
$K'=14_{27071}$) the $(-3,5,7)$--pretzel knot $P(-3,5,7)=15_{199038}$.
In fact, $P(-3,5,7)$ has also a diagram $D$ with $rb(D)>0$ (one comes
from a 32-crossing braid representation; see \cite{deg}). This shows,
that for some knots we may be lucky to find a diagram that excludes
sliceness directly, and working with $K'$ is not necessary.
Of course, by Fintushel-Stern's, and later Rudolph's work (see
\cite{Rudolph2}), the non-sliceness of this famous example has
been dealt with before. Note contrarily that, since both knots have
$\Dl=1$, they \em{are topologically} slice by Freedman's theorem.
\end{exam}

We found in total
5 prime knots up to 16 crossings (including the 3 so far mentioned),
for which the indirect Rudolph-Bennequin inequality proved essential
in excluding sliceness. All these knots have $\md_lP=2$, however, so
that on some more complicated diagram the direct inequality may apply~%
-- as seen for $P(-3,5,7)$. Note that, when the indirect argument
works, still \eqref{iy} holds, so that achiral knots cannot occur.
The attention to our examples was drawn by the problem of \cite{deg}
to find slice knots with $\md_lP>0$. Since the described method ruled
out all candidates for such a knot we had, the problem remains open.

\begin{remark}\label{r5.1}
Recently, Ozsvath and Szabo \cite{OS} defined a new ``signature''
invariant $\tau$ (for \em{knots}) using Floer homology, and Rasmussen
\cite{Rasmussen} a conjecturedly equivalent invariant $s$ using Khovanov
homology. This invariant lies between the two hand-sides of the 
slice Bennequin inequality $b(D)\le g_s(K)$. Thus it must confirm,
too, the non-sliceness of the examples in this section. However, it
is still non-trivial to calculate, and thus Rudolph-Bennequin numbers
remain a useful tool~-- in fact, one can estimate $s$ often easier
from them than calculating it directly. (Nonetheless Shumakovitch
\cite{Shu} computed $s$ on a number of knots with $\Dl=1$. He found
some knots of $s\ne 0$, where $\dl(P)=-4$ is relatively small, and
so the existence of non-trivial Rudolph-Bennequin numbers seems
unlikely.) On the opposite side, 
the calculation of $s$ is easy for alternating knots by virtue
of being equal to the usual signature $\sg$. This fact shows the
failure of the new invariant, too, for many of our examples,
including the first knot in figure \reference{fig1}, and the knots
in figures \reference{fig2} and \reference{fig4}.
It also explains the failure of the Bennequin numbers for such knots.
\end{remark}

\section{Fibered, prime, arborescent and hyperbolic knots\label{SFi}}

Motivated by Nakamura's construction \cite{Nakamura}, in \cite{aparb},
we gave another proof of Sakai-Kondo's result. (We found subsequently
that our construction was given, in a different context, also in
\cite{Murakami}.) The method explained in \cite{aparb} allows us to
choose our examples to have a few special properties. Above, we
considered a knot $K$, found by computation, and a knot $L$ of
unknotting number one with suitable Alexander polynomial. Then the
knots $K_k:=K\,\#\,\#^kL$ can be modified to prime knots $K_k'$
using tangle surgery. See proposition \reference{ppp} below.
Also, \cite{aparb} shows that one can choose $L$ fibered,
if its Alexander polynomial is monic. So one can obtain
(composite) fibered examples $K_k$ if $K$ is fibered. (How
to keep fiberedness in going from $K_k$ to $K_k'$ is not
clear at this point, though.)

\begin{prop}\label{ppp}
One can make a composite knot $K$ into a (prime)
hyperbolic knot $K'$ of arbitrarily large volume, preserving
$\Dl$, $H_1$, $\lm$, all non-singular Tristram-Levine signatures
$\sg_\xi$, and the Rasmussen invariant $s$. If $K$ is connected
sum of arborescent knots, one can also choose $K'$ to be arborescent.
\end{prop}

The same argument works for the Ozsvath-Szabo signature $\tau$
instead of $s$. If $\tau$ is equivalent to $s$, then there is
anyway nothing to do. If, however, $\tau$ is not equivalent to $s$,
and one would like to keep track of $\tau$ and $s$ \em{simultaneously},
an extra argument must be provided.

\proof In \cite{aparb} we showed that by tangle surgery
one can make $K$ into a hyperbolic knot of arbitrarily large
volume. If $K$ is arborescent connected sum, one can also
choose $K'$ to be arborescent. (More generally, the tangle surgery
can be chosen so as to preserve the Conway polyhedron of the diagram.)

Now, for tangle surgery one chooses the numbers of twists
to satisfy congruences modulo the determinant. Then $K$
and $K'$ have Goeritz matrices \cite{GL} $G$ and $G'$,
such that $d:=\dt(G)=\dt(G')$ and $G\equiv G'\bmod d$.
This implies that $H_1$ and $\lm$ are preserved.

The tangle surgeries preserve the Alexander polynomial.
It can also be easily observed that they preserve all
the non-singular Tristram-Levine signatures $\sg_\xi$.
Namely, by (in this or the reverse order) changing a
positive crossing (to become negative), applying
concordance, and changing a negative crossing, one 
obtains the original knot. So $\sg_\xi$ is changed by at most
$\pm 2$. But if $\Dl(\xi)\ne 0$, the sign of $\Dl(\xi)$
determines $\sg_\xi \bmod 4$; see \cite{tl}.
So the failure of $\sg_\xi$ to exhibit non-sliceness
persists under the tangle surgery even if the
Alexander polynomial has zeros on the unit circle.

Since $s$ is not connected to $\Dl$, so far tangle surgery
may alter it. We will argue how to remedy this problem.

Consider the pretzel knots of the form $P(-p,p+2,q)$, with $p,q>1$
odd and chosen so that $\Dl=1$. First note that all such knots have
$s=2$. Namely, by the main Theorem in \S 1 of \cite{Rudolph4},
these pretzel knots are quasipositive, and by proposition 5.3
of \cite{Rudolph4} have slice genus $1$. By \cite{Rasmussen}
one has then $s=2$. 

Now with a connected sum of a proper number $k$ of $P(-p,p+2,q)$
or their mirror images, we can make $s$ vanish. (We require that
$k$ is bounded by the number of surgeries in a way independent
on the number of twists in the surgered tangles.) In \cite{aparb} we
showed that the tangle surgery making the connected sum prime can
preserve the smooth concordance class (and hence $s$). 

To show hyperbolicity, now note that we can augment $p$ arbitrarily.
Similarly we can augment the dealternator twists in the surgered
tangles (also within a proper congruence class to keep track of
$H_1$ and $\lm$). Then by Thurston's hyperbolic surgery theorem
\cite{Thurston1} and the result of Adams on the hyperbolicity of
augmented alternating links \cite{Adams2}, the knots will be
hyperbolic for large number of crossings in the twists. To obtain
arbitrarily large volume, we augment the number of twist classes of
crossings, and use Adams' lower estimate on the volume of the augmented
alternating link in terms of the number of components \cite{Adams}.
(For an explanation, one may consult also \cite{Brittenham}.) \qed


{\bf Acknowledgement.}
Large parts of this paper were written during my stay at the
Advanced Mathematical Institute, Osaka City University. I would
like to thank Professor Akio Kawauchi for his invitation to
Osaka and the First KOOK Seminar International for Knot Theory
and Related Topics in Awaji-Shima (July 2004), and for all his 
support. Stefan Friedl, Tomomi Kawamura, Darren Long, Kunio
Murasugi and particularly Chuck Livingston contributed with
several useful remarks. However, the most substantial
improvements (incl. a complete reorganization of the introductory
part), were suggested by the referee. I owe to thank him and
editor Chuck Weibel for the care in dealing with my paper.

{\small
\newbox\@tempboxb
\setbox\@tempboxb=\hbox{P.\ Ozsv\'ath and Z.\ Szab\'o}

\let\old@bibitem\bibitem
\def\bibitem[#1]{\old@bibitem}


}

\end{document}